\documentclass[ECP,preprint]{ejpecp} 

\usepackage{xcolor}
\usepackage{mathrsfs}
\usepackage[inline, shortlabels]{enumitem}
 \usepackage{relsize}

\SHORTTITLE{A Mecke-type characterization of the Dirichlet--Ferguson measure}

\TITLE{A Mecke-type characterization\\of the Dirichlet--Ferguson measure\support{Research supported by the Sfb~1060 \emph{The Mathematics of Emergent Effects} (University of Bonn).
L.D.S.~gratefully acknowledges funding of his current position by the Austrian Science Fund (FWF) through project ESPRIT~208.
}\
    \thanks{The authors are grateful to Professor G\"unter Last and to two anonymous reviewers for some useful comments and suggestions.
    Results in this work were partly obtained during visits of the first named author to Swansea University.
    He is very grateful to the second named author and to the Department of Mathematics for their kind hospitality.}
    } 

\AUTHORS{%
  Lorenzo Dello Schiavo~\footnote{Institute of Science and Technology Austria, Austria.
    \EMAIL{lorenzo.delloschiavo@ist.ac.at}
    }
  \and 
Eugene Lytvynov\footnote{University of Swansea,
    United Kingdom. \EMAIL{e.lytvynov@swansea.ac.uk} }}


\AMSSUBJ{60G57} 

\SUBMITTED{} 
\ACCEPTED{} 

\VOLUME{0}
\YEAR{0}
\PAPERNUM{0}
\DOI{0}

\ABSTRACT{We prove a characterization of the Dirichlet--Ferguson measure over an arbitrary finite diffuse measure space.
We provide an interpretation of this characterization in analogy with the Mecke identity for Poisson point processes.}

\newtheorem*{theorem*}{Theorem}

\newcommand{\mfS}{{\mathfrak S}}

\newcommand{\mcD}{{\mathcal D}}
\newcommand{\mcG}{{\mathcal G}}
\newcommand{\mcM}{{\mathcal M}}

\newcommand{\mcP}{{\mathcal P}}
\newcommand{\mcQ}{{\mathcal Q}}

\newcommand{\mcX}{{\mathcal X}}

\newcommand{\msP}{{\mathscr P}}

\newcommand{\mbfM}{{\mathbf M}}
\newcommand{\mbfN}{{\mathbf N}}
\newcommand{\mbfX}{{\mathbf X}}

\newcommand{\mbfe}{{\mathbf e}}
\newcommand{\mbfs}{{\mathbf s}}
\newcommand{\mbfy}{{\mathbf y}}

\newcommand{\boldalpha}{{\boldsymbol \alpha}}

\newcommand{\mbbE}{{\mathbb E}}

\newcommand{\mbbX}{{\mathbb X}}

\newcommand{\emparg}{{\,\cdot\,}}

\newcommand{\Dir}[1]{D_{#1}}								

\DeclareMathOperator{\ev}{ev}

\DeclareMathOperator{\eqdef}{\colon\!\! \raisebox{-.4pt}{=}}

\newcommand{\rar}{\rightarrow}
\newcommand{\n}[1]{\overline{#1}}

\newcommand{\diff}{\mathop{}\!\mathrm{d}}						
	
\newcommand{\abs}[1]{\left\lvert#1\right\rvert}						
\newcommand{\set}[1]{\left\{#1\right\}}							
\newcommand{\tonde}[1]{\left(#1\right)}							
\newcommand{\ttonde}[1]{\big(#1\big)}							
\newcommand{\quadre}[1]{\left[#1\right]}							
\newcommand{\tquadre}[1]{\big[#1\big]}							

\newcommand{\seq}[1]{\tonde{#1}}								

\newcommand{\Mp}{{\mbfM}}

\newcommand{\pfwd}{*}

\DeclareMathOperator{\N}{{\mathbb N}}
\DeclareMathOperator{\R}{{\mathbb R}}

\DeclareMathOperator{\law}{law}

\newcommand{\comma}{\;\textrm{,}\quad}

\newcommand{\fstop}{\;\textrm{.}}

\newcommand{\compo}{\diamond}

\let\epsilon\varepsilon
\let\subset\subseteq

\let\temp\phi
\let\phi\varphi
\let\varphi\temp

\newcommand{\Beta}{\mathrm{B}}

\newcommand{\DF}{{\mcD}}
\newcommand{\GP}{{\mcG}}
\newcommand{\PP}{{\mcP}}
\newcommand{\Gam}[1]{G_{#1}}
\newcommand{\BetaD}[1]{B_{#1}}

\newcommand{\sym}[1]{{\scriptscriptstyle{(#1)}}}

\allowdisplaybreaks

\begin{document}


\section{Introduction and the main result}
Let~$(\mbbX,\mcX)$ be a measurable space with $\sigma$-algebra~$\mcX$, and assume that points in~$\mbbX$ are measurable, i.e.~$\set{x}\in\mcX$ for every~$x\in\mbbX$.
For any (non-negative) measure~$\mu$ on~$(\mbbX,\mcX)$, and any $\mcX$-measurable~$f\colon \mbbX\to \overline{\R}\eqdef\R\cup\set{\pm\infty}$ we denote by~$\mu f$ the integral of~$f$ with respect to~$\mu$, whenever this makes sense.


We denote by~$\Mp\eqdef \Mp(\mbbX)$ the cone of all measures on~$\mbbX$ with values in~$[0,+\infty]$, and by~$\R_b(\mbbX)$, respectively~$\R_+(\mbbX)$ the linear space of all bounded, respectively non-negative, $\mcX$-measurable $\R$-valued functions.
We endow~$\Mp$ with the coarsest $\sigma$-algebra~$\mcM$ for which all functions of the form~$\mu\mapsto \mu B\in [0,+\infty]$ with~$B\in\mcX$ are measurable for the Borel $\sigma$-algebra of the extended half-line.
We denote by~$\mbfN\eqdef \mbfN(\mbbX)$ the space of all $\overline{\N}_0$-valued elements of~$\mbfM$. Here and elsewhere we set~$\overline\N_0\eqdef \N_0\cup \set{+\infty}$. It holds that~$\mbfN \in\mcM$.

When~$(\mbbX,\mcX)$ is a Borel space, for each~$\gamma\in\mbfN$ there exists an at most countable family of not necessarily distinct points~$x_i\in \mbbX$ such that~$\gamma=\sum_{i=1}^{ \gamma \mbbX} \delta_{x_i}$, see e.g.~\cite[Prop.~6.2]{LasPen17+}.

By a~\emph{random measure} on~$\mbbX$ we mean any~$(\Mp,\mcM)$-valued random field.
If a random measure~$\gamma$ is concentrated on~$\mbfN$, we say that it is a~\emph{point process} in~$\mbbX$.
Following~\cite[\S2.1]{LasPen17+} we say that a measure in~$\mbfM$ is \emph{$s$-finite} if it is the sum of at most countably many \emph{finite} measures in~$\mbfM$.
Everywhere in the following let~$\sigma$ be an $s$-finite element of~$\Mp$.
A random measure~$\nu$ has \emph{intensity} (\emph{measure})~$\sigma$ if
\begin{equation}\label{eq:Intensity}
\mbbE\tquadre{\nu B}= \sigma B \comma \qquad B\in\mcX\fstop
\end{equation}

\paragraph{Poisson measures}
Among all point processes, a remarkable and ubiquitous example is given by the \emph{Poisson point process}~$\gamma$ with intensity~$\sigma$, i.e., the point process in~$\mbbX$ with Laplace transform
\begin{align}\label{eq:LapPP}
\mbbE\quadre{ e^{-\gamma f} }= \exp\quadre{-\int \tonde{1-e^{f(x)}}  \diff\sigma(x)}, \qquad  f\in \R_+(\mbbX) \comma
\end{align}
see e.g.~\cite[\S3.1]{Kin93} or \cite[Thm.~3.9]{LasPen17+}.

We denote by~$\PP_\sigma$ the law of a Poisson point process with intensity~$\sigma$, and we write~$\gamma\sim\PP_\sigma$ to indicate that~$\gamma$ is distributed as~$\PP_\sigma$.
Recall the following characterization of~$\gamma\sim \PP_\sigma$, usually known as the \emph{Mecke identity}.
\begin{theorem*}[Mecke identity for~$\PP_\sigma$ {\cite[Satz~3.1]{Mec67},~\cite[Thm.~4.1]{LasPen17+}}]
Let~$\mbbX$ and $\sigma$ be as above and let $\gamma$ be a random measure over~$\mbbX$. Then, the following statements are equivalent:

\begin{enumerate}
\item[$(i)$] $\gamma$ is a Poisson point process with intensity~$\sigma$;

\item[$(ii)$]  for every measurable function~$F\colon \Mp\times \mbbX\rar [0,+\infty)$,
\begin{align}\label{eq:MeckePP}
\mbbE\quadre{ \int F(\gamma, x)  \diff\gamma(x)} = \int \mbbE\tquadre{F(\gamma+\delta_x, x)}  \diff \sigma(x) \fstop
\end{align}
\end{enumerate}
\end{theorem*}


\paragraph{Main result}
Let~$\msP$ be the subset of \emph{probability measures} in~$\Mp$.
It holds that $\msP\in\mcM$.
If a random measure is almost surely an element of~$\msP$, we say that it is a~\emph{random probability} (\emph{measure}).

The aim of this work is to show how the law~$\DF_\sigma$ of a Dirichlet--Ferguson process (see~\S\ref{s:Prel} below) may be regarded as the natural analog of the Poisson measure~$\PP_\sigma$ when one replaces~$\mbfN$ with~$\msP$.

\begin{theorem}[A Mecke-type characterization of~$\DF_\sigma$]\label{t:MeckeDF}
Let $\sigma$ be a \emph{finite} \emph{diffuse} (i.e., atomless) measure on~$(\mbbX,\mcX)$, and set~$\beta\eqdef \sigma \mbbX$. Then, for any random measure~$\eta$ over~$\mbbX$, the following statements are equivalent:
\begin{enumerate}[$(i)$]
\item \label{i:1} $\eta$ is a Dirichlet--Ferguson process on~$\mbbX$ with intensity~$\sigma$;
\item \label{i:2} for every measurable function~$G\colon \Mp\rar \R_+$,
\begin{align}\label{eq:MeckeDF}
\mbbE\tquadre{ \eta\mbbX\, G(\eta)} = \int \int_0^1 \mbbE \quadre{G\big((1-t)\eta +t\delta_x\big)} (1-t)^{\beta-1}\diff t \diff\sigma(x) \fstop
\end{align}
\end{enumerate}

Moreover, if~$\eta$ is a Dirichlet--Ferguson process on~$\mbbX$ with intensity~$\sigma$, then for every non-negative (or bounded) measurable function $F\colon\msP\times \mbbX\rar \R$,
\begin{align}\label{vcghfsky}
\mbbE\quadre{ \int F(\eta, x)  \diff\eta(x) }=\int \int_0^1\mbbE\quadre{F\big((1-t)\eta+t\delta_x, x\big)} (1-t)^{\beta-1}\diff t \diff\sigma(x) \comma
\end{align}
and  for every non-negative (or bounded) measurable function $R\colon\msP\times \mbbX\times [0,1]\rar \R$,
\begin{align}\label{eq:MeckeDFCor}
\mbbE\quadre{ \int R\big(\eta, x, \eta(x)\big) \diff\eta(x)}= \int \int_0^1 \mbbE\quadre{R\big((1-t)\eta+t\delta_x, x, t\big)} (1-t)^{\beta-1}\diff t \diff\sigma(x) \fstop
\end{align}
In formula \eqref{eq:MeckeDFCor} we set $\eta(x)\eqdef\eta(\{x\})$.
\end{theorem}

\begin{remark}
Formula \eqref{eq:MeckeDF} is a special case of formula~\eqref{vcghfsky}, while formulas~\eqref{vcghfsky} and~\eqref{eq:MeckeDFCor} are in fact equivalent.
Formulas~\eqref{vcghfsky} and~\eqref{eq:MeckeDFCor} are an integral reformulation of the celebrated \emph{stick-breaking} construction of~$\DF_\sigma$ obtained by J.~Sethuraman,~\cite[\S3]{Set94}.
\end{remark}

\begin{remark}
By Theorem~\ref{t:MeckeDF}, the law of random probability measure~$\eta$ on~$\mbbX$ is the Dirichlet--Ferguson measure~$\DF_\sigma$ with intensity~$\sigma$ if and only if, for every measurable function~$G\colon \msP\rar \R_+$
\begin{align*}
\mbbE\tquadre{G(\eta)} = \int \int_0^1 \mbbE \quadre{G\big((1-t)\eta +t\delta_x\big)} (1-t)^{\beta-1}\diff t \diff\sigma(x) \fstop
\end{align*}
\end{remark}

\begin{remark}Note that, when $\sigma$ is a probability measure on $\mbbX$, in formulas \eqref{eq:MeckeDF}--\eqref{eq:MeckeDFCor}, the factor $(1-t)^{\beta-1}$ becomes 1.
If~$\sigma$ is not a probability measure, let $\n\sigma\eqdef \sigma/\beta\in \msP$, denote by~$\Beta$ the Beta function and by
\[
\diff\BetaD{a,b}(t)\eqdef \frac{t^{a-1} (1-t)^{b-1} \diff t}{\Beta(a,b)}
\]
 the Beta distribution on~$[0,1]$ with shape parameters $a>0$ and $b>0$.
Then we have the following equality of the probability measures on $\mbbX\times[0,1]$:
\[
\diff \sigma(x)\, (1-t)^{\beta-1} \diff t
=\diff \n\sigma(x)\, \beta (1-t)^{\beta-1} \diff t
=\diff\n \sigma(x)\,\diff\BetaD{1,\beta}(t) \fstop
\]
\end{remark}

\begin{remark}\label{r:Preprint}
A first version of this paper appeared as the arXiv preprint 1706.07602. The main difference with the present work is that the space~$\mbbX$ in the preprint was assumed to be a locally compact Polish space.
After the publication of~arXiv:1706.07602, G.~Last~\cite{Las19} proved a characterization of the Dirichlet--Ferguson process that uses a weaker form of equality~\eqref{vcghfsky}.
More precisely, by~\cite[Thm.~1.5]{Las19}, if a random measure satisfies, for all measurable functions~$F\colon\mbfM\times\mbbX\to [0,+\infty]$,
\[
\mbbE\quadre{\int F(\eta, x)\diff\eta(x)}=\int \int_0^1 \mbbE\quadre{F\ttonde{(1-t)\eta+t\delta_x,x}}\diff G(t) \diff\n\sigma(x)
\]
for some probability measure~$G$ on~$[0,1]$, then~$G=\BetaD{1,\beta}$ and so~$\eta$ is a Dirichlet--Ferguson process with intensity~$\sigma$.
Note that this result does not imply Theorem~\ref{t:MeckeDF} since formula~\eqref{eq:MeckeDF} is a special case of~\eqref{vcghfsky}.
\end{remark}

\begin{remark}
Other characterizations of the Dirichlet--Ferguson measure and of the Dirichlet distribution are also available, e.g.:~\cite{RegGugDiN02} for characterization of (functionals) of~$\DF_\sigma$ via contour-integral methods,~\cite{JiaDicKuo04} for a characterization via $c$-transform, and~\cite{LzDS19} for a characterization via Fourier transform.
\end{remark}



\paragraph{Motivations and applications}
The Mecke identity~\eqref{eq:MeckePP} and  its generalization to other random measures, e.g.\ the \emph{Georgii--Nguyen--Zessin formula} for Gibbs measures~\cite{Geo76, NguZes79, MatWarMec78}, have important applications in the theory of point processes and stochastic dynamics of interacting particle systems, see e.g.~\cite{dSKR,MR}. 
In a similar fashion, when~$\eta$ is a measure-valued L\'evy process, suitable Mecke-type identities for~$\eta$ have been a key tool in the study of stochastic dynamics of measure-valued diffusion processes having the law of~$\eta$ as their invariant measure.
Among such Mecke-type identities we have also
identity~\eqref{eq:MeckeGP} for the gamma measure~$\GP_\sigma$ (see~\S\ref{s:Gamma} below), used to establish an integration-by-parts formula for the $\GP_\sigma$-invariant diffusion,~\cite{CKL,KonLytVer15}.

After the publication of this work as the arXiv preprint 1706.07602, several applications of our main result have appeared. We briefly summarize them in~\S\ref{s:Applications}.

\paragraph{Plan of the work}
Below, in \S\ref{s:Prel}, we discuss preliminary notions and facts, and in  \S\ref{vufuddx} we prove Theorem~\ref{t:MeckeDF} and discuss several corollaries, including a characterization of the Dirichlet distribution.

\section{Preliminaries}\label{s:Prel} 

\paragraph{The Dirichlet--Ferguson measure}
For integer $k\ge 2$, let $\Delta^{k-1}$ denote the  standard closed $(k-1)$-dimensional simplex in~$\R^{k}$, i.e.,
\[
\Delta^{k-1}\eqdef \big\{(y_1,\dotsc,y_{k}): y_i\ge 0,\ y_1+\dots+y_{k}=1\big\}\fstop
\]
Write $\R_+\eqdef (0,\infty)$.
For $\boldalpha\eqdef \seq{\alpha_1,\dotsc, \alpha_{k}}$ in~$\R^{k}_+$\,,  the \emph{Dirichlet distribution} with parameter~$\boldalpha$ is the probability measure  on $\Delta^{k-1}$ denoted by $\Dir{\boldalpha}$ and defined by
\begin{equation}\label{bvyuf7r8o}
\Dir{\boldalpha}(A)\eqdef  \frac{1}{\Beta(\boldalpha)}\int_A \left(\prod_{i=1}^{k} y_i^{\alpha_i-1}\right) \diff\mathcal H^{k-1}(y_1,\dotsc, y_k)
\end{equation}
for each measurable subset $A$ of $\Delta^{k-1}$, where~$\mathcal H^{k-1}$ denotes the Hausdorff measure on~$\Delta^{k-1}$. In formula \eqref{bvyuf7r8o},  $\Beta(\cdot)$ is the multivariate Beta function.

For integer $k\ge2$, an \emph{ordered partition} is a vector~$\seq{X_1,\dotsc, X_k}$ with the following properties:
\begin{enumerate}[$(a)$]
\item $X_i\in\mcX$ is a non-empty measurable subset of~$\mbbX$ for each $i\leq k$;
\item the sets~$X_i$ form a partition of~$\mbbX$, i.e.\ they are pairwise disjoint and their union coincides with~$X$.
\end{enumerate}
We denote by $\mathfrak P_{k}(\mbbX)$ the set of ordered partitions $\mbfX\eqdef\seq{X_1,\dotsc, X_{k}}$ of $\mbbX$.
For each $\mbfX\in \mathfrak P_{k}(\mbbX)$, we  define the evaluation map~$\ev_\mbfX \colon \mbfM\to\R_+^k$ by
\[
\ev_\mbfX\colon \eta \longmapsto \seq{\eta X_1,\dotsc, \eta X_{k}} \fstop
\]
Note that the map~$\ev_\mbfX$  is measurable, and so is its restriction~$\ev_\mbfX \colon \msP \to\Delta^{k-1}$.

\begin{lemma}\label{l:Identification}
Let~$\mcP$ and~$\mcQ$ be probability measures on~$(\mbfM,\mcM)$. Then, $\mcP=\mcQ$ if and only if~$(\ev_\mbfX)_\pfwd \mcP=(\ev_\mbfX)_\pfwd \mcQ$ for every~$\mbfX\in\mathfrak{P}_k(\mbbX)$ and every~$k\geq 2$.
\begin{proof}
The forward implication is trivial, thus it suffices to show the reverse one.
Since~$\mcP$ and~$\mcQ$ are probability measures, by e.g.~\cite[Lem.~I.9.4]{Bog07}, it suffices to verify that~$\mcP$ and~$\mcQ$ coincide on any algebra of sets generating $\sigma$-algebra~$\mcM$.
By definition of push-forward measure, for each Borel measurable~$A\subset \R_+^k$,
\begin{align*}
\mcP\quadre{\ev_\mbfX^{-1}(A)}=\mcQ\quadre{\ev_\mbfX^{-1}(A)} \comma
\end{align*}
and the conclusion follows since the~$\mcM$ is generated by all sets of the form~$\ev_\mbfX^{-1}(A)$ with~$A\subset \R_+^k$ Borel measurable, $k\ge2$, and~$\mbfX\in\mathfrak{P}_k(\mbbX)$.
\end{proof}
\end{lemma}

The \emph{Dirichlet--Ferguson process} $\eta$ with intensity~$\sigma$, see~\cite{Fer73}, is the unique random probability measure over~$\mbbX$ satisfying the following two conditions:
\begin{enumerate}[$(a)$]

\item for each $B\in\mcX$ with $\sigma B=0$, we have $\eta B=0$ a.s.;

\item  for each integer $k\ge2$ and $\mbfX\in \mathfrak P_{k}(\mbbX)$ additionally so that~$\sigma X_i>0$ for each~$i\leq k$,
\begin{align}\label{bur8ut}
\ev_\mbfX(\eta) \sim \Dir{\ev_\mbfX(\sigma)} \comma
\end{align}
i.e., the push-forward of  $\DF_\sigma$ under $\ev_\mbfX$ is equal to the Dirichlet distribution~$\Dir{\ev_\mbfX(\sigma)}$ with parameter~$(\sigma X_1,\dotsc, \sigma X_k)$.
\end{enumerate}

\paragraph{The gamma measure}\label{s:Gamma}
A compound Poisson process~$\nu$ on~$\mbbX$ is a random measure on~$(\mbbX,\mcX)$ of the form
\[
\nu A = \int_{\R_+} \int_A s \diff\gamma(x,s) \comma
\]
where~$\gamma$ is a Poisson point process on~$\hat\mbbX\eqdef\mbbX\times \R_+$ with a given intensity~$\hat\sigma$.
When
\begin{align}\label{eq:IntensityGamma}
\hat\sigma=\sigma\otimes s^{-1}e^{-s}\diff s\comma
\end{align}
the corresponding compound Poisson process is termed \emph{gamma process} with intensity~$\sigma$, and we denote its law by~$\GP_\sigma$.
It has Laplace transform
\[
\mbbE\quadre{e^{-\nu f}} = \exp\quadre{-\int \log\big(1+f(x)\big) \diff\sigma(x)}\comma\qquad  f\in \R_+(\mbbX) \comma
\]
see \cite[Eqn.~(7)]{TsiVerYor01} or~\cite[Example~15.6]{LasPen17+}.

\begin{lemma}[Mecke identity for the gamma measure]\label{cdtrde4ew54}
Let~$\nu$ be a gamma process on~$\mbbX$ with \emph{finite diffuse} intensity~$\sigma$. 
%
Then, for every measurable function $G\colon \Mp\times \mbbX\rar \R_+$,
\begin{align}\label{eq:MeckeGP}
\mbbE\quadre{ \int G(\nu, x) \diff\nu(x)}= \int \int_0^\infty \mbbE\tquadre{G(\nu+s\delta_x,x)} e^{-s}\diff s\diff\sigma(x) \fstop
\end{align}

\begin{proof}
Let~$\hat\sigma$ be as in~\eqref{eq:IntensityGamma} and~$\gamma\sim\PP_{\hat\sigma}$ be a Poisson point process on~$\mbbX\times\R_+$ with intensity~$\hat\sigma$.
By definition of the gamma process and by~\eqref{eq:MeckePP} for~$\gamma$, for every measurable function~$G\colon \mbfM\times \mbbX\to \R_+$,
\begin{align*}
\mbbE&\quadre{\int G(\nu,x)\diff\nu(x)}=
\\
&= \mbbE\quadre{\int_{\mbbX\times\R_+} G\tonde{\int_0^\infty s'\diff\gamma(\emparg, s'),x} s \diff\gamma(x,s)}
\\
&= \mbbE\quadre{\int_{\mbbX\times\R_+}  G\tonde{\int_0^\infty s'\diff\gamma(\emparg, s')+\int_{\R_+} s' \diff\delta_{(x,s)}(\emparg,s'),x} s \diff\hat\sigma(x,s)}
\\
&=\int_\mbbX \int_0^\infty \mbbE\quadre{G(\nu+s\delta_x,x)} e^{-s}\diff s \diff\sigma(x) \fstop \qedhere
\end{align*}
\end{proof}
\end{lemma}

It was shown in~\cite[\S4, Thm.~2, p.~219]{Fer73} (see also \cite[Lem.~1]{TsiVer99}) that the Dirichlet--Ferguson measure~$\DF_\sigma$ is the `simplicial part' of~$\GP_\sigma$.
More precisely, denote by  $\Gamma(\cdot)$  the  gamma function and by
\[
\diff\Gam{k,\theta}(s)\eqdef \frac{\theta^{-k}}{\Gamma(k)} s^{k-1} e^{-\tfrac{s}{\theta}}\diff s
\]
the gamma distribution on~$\R_+$ with shape parameter~$k\in\R_+$ and scale parameter~$\theta\in\R_+$.
Further consider the bi-measurable bijection $\mathfrak R:\Mp\setminus\set{0}\to\msP\times\R_+$ given by
\[
\mathfrak R(\nu)\eqdef\left(\frac{\nu}{\nu \mbbX},\,\nu \mbbX\right)\fstop
\]
Then,~$\nu$ is a gamma process with intensity~$\sigma$ if and only if~$\nu/\nu\mbbX$ is a Dirichlet--Ferguson process with intensity~$\sigma$, the random fields~$\nu/\nu\mbbX$ and~$\nu\mbbX$ are independent, and~$\nu\mbbX$ is a $G_{\beta,1}$-distributed $\R_+$-valued random variable.
(Recall that~$\beta=\sigma\mbbX$.)
Equivalently,
\begin{equation}\label{eq:Simplicial}
\mathfrak R_\pfwd\, \GP_\sigma=\DF_\sigma\otimes\Gam{\beta,1},
\end{equation}
i.e., the push-forward of $\GP_\sigma$ under $\mathfrak R$ is the product measure $\DF_\sigma\otimes\Gam{\beta,1}$.




\section{Proof and  corollaries}\label{vufuddx}

\begin{proof}[Proof of Theorem~\ref{t:MeckeDF}]
We first prove that a Dirichlet--Ferguson process~$\eta$ satisfies formula \eqref{vcghfsky}, hence \eqref{eq:MeckeDF}.
Let~$\eta\sim\DF_\sigma$, and~$R$ be a random variable with values in~$\R_+$ and gamma distribution~$\Gam{\beta,1}$, independent of the random probability~$\eta$.
Note that~$\mbbE[R]=\frac{\Gamma(\beta+1)}{\Gamma(\beta)}$.
We define a random measure~$\nu\eqdef R\eta$.
By~\eqref{eq:Simplicial} we have~$\nu\sim \GP_\sigma$ and, by construction,~$R=\nu\mbbX$ and~$\eta=\frac{\nu}{\nu\mbbX}$.
Then we have
\begin{align*}
\mbbE\quadre{\int F(\eta, x ) \diff\eta(x)}=&\ \frac{\Gamma(\beta)}{\Gamma(\beta+1)} \mbbE\quadre{R \int F(\eta, x) \diff\eta(x)}
\\
=&\ \frac{\Gamma(\beta)}{\Gamma(\beta+1)} \mbbE\quadre{\int F(\eta, x)\, R \diff\eta(x)}
\\
=&\ \frac{\Gamma(\beta)}{\Gamma(\beta+1)} \mbbE\quadre{\int F\tonde{\frac{\nu}{\nu\mbbX}, x}\,  \diff\nu(x)} \fstop
\end{align*}
Using Lemma~\ref{cdtrde4ew54}, we continue the above chain of equalities as follows:
\begin{align*}
&=\frac{\Gamma(\beta)}{\Gamma(\beta+1)}  \mbbE\quadre{\int \int_0^\infty F\left(\frac{\nu+s\delta_x}{\nu \mbbX+s}\,,x\right) e^{-s}\diff s \diff\sigma(x)}
\\
&=\frac{\Gamma(\beta)}{\Gamma(\beta+1)} \mbbE\quadre{\int \int_0^\infty F\left(\frac{R\eta+s\delta_x}{R+s}\,,x \right) e^{-s}\diff s \diff\sigma(x)}
\\
&=\frac{\Gamma(\beta)}{\Gamma(\beta+1)} \int \int_0^\infty  \int_0^\infty \mbbE\quadre{F\left(\frac{r\eta+s\delta_x}{r+s}\,,x \right)} e^{-s}\diff s \diff\Gam{\beta,1}(r)\diff\sigma(x)
\\
&=\frac{1}{\Gamma(\beta+1)} \int \int_0^\infty \int_0^\infty \mbbE\quadre{F\left(\frac{r}{r+s}\,\eta+\frac{s}{r+s}\,\delta_x,x \right)} e^{-s}\diff s\ r^{\beta-1} e^{-r} \diff r \diff\sigma(x)\comma
 \intertext{whence the change of variable $t=\tfrac{s}{r+s}$ (for a fixed $s$) yields}
 &=\frac{1}{\Gamma(\beta+1)} \int \int_0^\infty \int_0^1\, s^{\beta-1}\frac{(1-t)^{\beta-1}}{t^{\beta-1}} e^{-\frac{s(1-t)}{t}}\, \mbbE\quadre{F\big((1-t)\eta+t\delta_x,x\big)} \frac{ s \, \diff t}{t^2} e^{-s}\diff s \diff\sigma(x)
 \\
 &=\frac{1}{\Gamma(\beta+1)} \int \int_0^1 \frac{ (1-t)^{\beta-1}}{t^{\beta+1}} \,  \mbbE\quadre{F\big((1-t)\eta+t\delta_x,x\big)} \diff\sigma(x) \int_0^\infty e^{-s}s^\beta e^{-\frac{s(1-t)}{t}} \diff s \diff t\\
&=\int \int_0^1 \mbbE\tquadre{F\big((1-t)\eta+t\delta_x,x\big)} (1-t)^{\beta-1}\diff t \diff\sigma(x)\fstop
\end{align*}
To prove formula \eqref{eq:MeckeDFCor}, choose $F(\eta, x)= R\big(\eta, x, \eta(x)\big)$ in \eqref{vcghfsky}, which gives
\begin{align*}
&\mbbE \quadre{ \int R\big(\eta, x, \eta(x)\big) \diff\eta(x)}=
\\
&= \int \int_0^1 \mbbE\quadre{R\big((1-t)\eta+t\delta_x,x,(1-t)\eta(x)+t\big)}  (1-t)^{\beta-1}\diff t \diff\sigma(x)
\\
&=\int \int_0^1\mbbE\quadre{R\big((1-t)\eta+t\delta_x,x,t\big)} (1-t)^{\beta-1}\diff t \diff\sigma(x)\comma
\end{align*}
where we used that, for a fixed $\eta\in\msP$, we have $\eta(x)=0$ for $\sigma$-a.e.\ $x\in \mbbX$, as a consequence of~\eqref{eq:Intensity} for~$\eta\sim\DF_\sigma$ in place of~$\nu\sim\GP_\sigma$.

\medskip

For the reverse implication we consider a random measure~$\eta$ over~$\mbbX$ that satisfies~\eqref{eq:MeckeDF}.
We need to show that~$\law\eta=\DF_\sigma$.
Let us first show that~$\eta\in\msP$ a.s. Choosing $G\equiv 1$ in~\eqref{eq:MeckeDF}, we get
\begin{equation}\label{cyteri67}
\mbbE\tquadre{\eta \mbbX}=1 \fstop
\end{equation}
In particular, $\eta \mbbX<\infty$ a.s.
Next, choosing $G(\eta)=\eta \mbbX$ in~\eqref{eq:MeckeDF} and using \eqref{cyteri67}, we get
\begin{align}\label{dre6u}
\mbbE\quadre{(\eta \mbbX)^2}=&\int \int_0^1\mbbE\tquadre{(1-t)\,\eta \mbbX +t} (1-t)^{\beta-1} \diff t\diff\sigma(x)
\\
=& \int \diff\sigma(x) \int_0^1\big((1-t)+t\big)(1-t)^{\beta-1}  \diff t  =1\fstop\notag
\end{align}
By \eqref{cyteri67} and \eqref{dre6u}, the random variable $\eta \mbbX$ has zero variance, hence it is deterministic.
Thus, $\eta \mbbX=1$ a.s., so $\eta\in\msP$ a.s.
Hence, formula \eqref{eq:MeckeDF} becomes
 \begin{align}\label{fvgyugwe}
\mbbE\tquadre{G(\eta)} = \int\int_0^1 \mbbE\quadre{G\big((1-t)\eta +t\delta_x\big) } (1-t)^{\beta-1} \diff t \diff\sigma(x)\comma
\end{align}
and it holds for every measurable bounded function $G:\msP \to\R$.

Let $B\in\mcX$ be such that $\sigma B=0$. By \eqref{fvgyugwe},
\begin{align*}
\mbbE\quadre{\eta B}=&\int\int_0^1 \mbbE\tquadre{(1-t)\, \eta B +t\mathbf 1_B(x)} (1-t)^{\beta-1}\diff t \diff\sigma(x)
\\
=&\ \mbbE\quadre{\eta B} \beta \int_0^1 (1-t)^{\beta} \diff t+\sigma B\int_0^1 t (1-t)^{\beta-1} \diff t
\\
=&\ \frac{\beta}{\beta+1}\mbbE\tquadre{\eta B}+0\comma
\end{align*}
which implies
\[
\mbbE\tquadre{\eta B} =0\fstop
\]
Hence, $\eta B =0$ a.s.

Let $k\ge2$ and  $\mbfX=(X_1,\dots,X_{k})\in \mathfrak P_{k}(\mbbX)$ be an ordered partition with~$\sigma X_i>0$ for all~$i\leq k$.
In order to prove that $\eta\sim\DF_\sigma$, it remains to show that the distribution of the random vector $\ev_\mbfX(\eta)$ in $\R^{k}$ (in fact, in $\Delta^{k-1}$) is~$\Dir{\ev_\mbfX(\sigma)}$.

We recall that the \emph{Hadamard product} $\diamond:\R^k\times\R^k\to\R^k$ is defined by
\[
\mbfs^\sym{1}\compo\mbfs^\sym{2}\eqdef(s^\sym{1}_1s^\sym{2}_1,\dots,s^\sym{1}_ks^\sym{2}_k),\qquad \mbfs^\sym{i}=(s^\sym{i}_1,\dots,s^\sym{i}_k)\in\R^k,\quad i=1,2 \fstop
\]
This binary operation is obviously associative and commutative.


Write~$\boldalpha\eqdef\ev_\mbfX(\sigma)$.  Fix any $\mathbf s=(s_1,\dots,s_{k})\in\R^{k}$, and let $g(x)\eqdef\sum_{i=1}^{k}s_i\mathbf 1_{X_i}(x)\in \R_b(\mbbX)$. 
Then
\begin{equation}\label{fst5w}
\eta g=\mathbf s \cdot \ev_\mbfX(\eta),\qquad \eta\in \msP\comma
\end{equation}
and
\begin{equation}\label{fvte6e54w}
\sigma(g^n)=\mbfs^{\compo n}\cdot \boldalpha,\qquad n\in\N_0 \fstop
\end{equation}
For $n\in\mathbb N$, we get, by \eqref{fvgyugwe}--\eqref{fvte6e54w},
\begin{align*}
\mbbE\quadre{\big(\mathbf s \cdot \ev_\mbfX(\eta)\big)^n}=& \int \int_0^1 \mbbE\quadre{\big((1-t)\, \eta g+tg(x)\big)^n} (1-t)^{\beta-1}\diff t \diff\sigma(x)
\\
=&\sum_{i=0}^n\binom{n}{i} \mbbE\quadre{(\eta g)^i} \int g(x)^{n-i}  \diff\sigma(x) \int_0^1 (1-t)^{\beta+i-1} t^{n-i}\diff t
\\
=&\sum_{i=0}^n\binom ni 
\Beta(\beta+i,n-i+1)\, \mbbE\quadre{\big(\mathbf s \cdot \ev_\mbfX(\eta)\big)^i (\mbfs^{\compo (n-i)}\cdot \boldalpha)}
\\
=& \sum_{i=0}^n \frac{n! \, \Gamma(\beta+i)}{i!\, \Gamma(\beta+n+1)}\, \mbbE\quadre{\big(\mathbf s \cdot \ev_\mbfX(\eta)\big)^i (\mbfs^{\compo (n-i)}\cdot \boldalpha)}
\\
=&\sum_{i=0}^n \frac{(n)_{n-i}}{(\beta+n)_{n+1-i}}\, \mbbE\quadre{\big(\mathbf s \cdot \ev_\mbfX(\eta)\big)^i (\mbfs^{\compo (n-i)}\cdot \boldalpha)}
\\
=&\frac\beta{\beta+n}\, \mbbE\quadre{\big(\mathbf s \cdot \ev_\mbfX(\eta)\big)^n}
\\
&+\sum_{i=0}^{n-1} \frac{(n)_{n-i}}{(\beta+n)_{n+1-i}}\, \mbbE\quadre{\big(\mathbf s \cdot \ev_\mbfX(\eta)\big)^i} (\mbfs^{\compo (n-i)}\cdot \boldalpha) \comma
\end{align*}
where $(r)_k$ denotes the falling factorial: $(r)_0\eqdef1$ and $(r)_k\eqdef r(r-1)\dotsm(r-k+1)$ for~$k\in\mathbb N$.
Therefore,
\begin{align}
\label{bgufr7rd}
\mbbE\quadre{ \big(\mathbf s \cdot \ev_\mbfX(\eta)\big)^n} & =
\frac{\beta+n}{n}\sum_{i=0}^{n-1} \frac{(n)_{n-i}}{(\beta+n)_{n+1-i}}\, \mbbE\quadre{\big(\mathbf s \cdot \ev_\mbfX(\eta)\big)^i} (\mbfs^{\compo (n-i)}\cdot \boldalpha)
\\
& =
\sum_{i=0}^{n-1} \frac{(n-1)_{n-1-i}}{(\beta+n-1)_{n-i}}\, \mbbE\quadre{\big(\mathbf s \cdot \ev_\mbfX(\eta)\big)^i} (\mbfs^{\compo (n-i)}\cdot \boldalpha) \fstop\notag
\end{align}
The recurrence relation \eqref{bgufr7rd}, with initial condition~$\mbbE\tquadre{\ttonde{\mathbf s \cdot \ev_\mbfX(\eta)}^0}=1$ for every~$\mbfs$ and~$\mbfX$, uniquely determines the moments 
\begin{equation}\label{vfr7}
\mbbE\quadre{\big(\mathbf s \cdot \ev_\mbfX(\eta)\big)^n} \comma \qquad \mathbf s\in\R^{k}\comma n\in\mathbb N_0 \fstop\end{equation}

By the polarization identity (see e.g.~\cite[Ch.~2\, \S2.1, Eqn.~(2.7)]{BerKon95}, p.~132), we therefore conclude that the recurrence relation~\eqref{bgufr7rd} uniquely determines the moments
\[
\mbbE\quadre{\ttonde{\mathbf s_1 \cdot \ev_\mbfX(\eta)}\cdots \ttonde{\mathbf s_n \cdot \ev_\mbfX(\eta)}} \comma \qquad \mathbf s_1,\dotsc, \mathbf s_n\in\R^{k}\comma n\in\mathbb N_0 \fstop
\]
In turn, this implies that the moments
\begin{equation}\label{eq:Moments2}
\mbbE\quadre{(\eta X_1)^{j_1}\cdots (\eta X_k)^{j_k}}\comma \qquad j_1,\dotsc, j_k\in \N_0\comma
\end{equation}
too are uniquely determined by~\eqref{bgufr7rd}.
Since~$\eta$ is a random probability measure,~$(\ev_\mbfX)_\pfwd \law \eta$ is supported on the unit simplex~$\Delta^{k-1}$.
Thus, each of the moments in~\eqref{eq:Moments2} is bounded in modulus by~$1$.
Hence, by~\cite[Ch.~8, \S5.7, Thm.~5.10, p.~735]{Ber} there exists a unique measure~$\mu_\mbfX$ on~$\R^k$ such that
\[
\int_{\R^k} t_1^{j_1}\cdots t_k^{j_k} \diff\mu_\mbfX(t_1,\dotsc, t_k) = \mbbE\quadre{(\eta X_1)^{j_1}\cdots (\eta X_k)^{j_k}}\fstop
\]
Since the Dirichlet--Ferguson process satisfies~\eqref{eq:MeckeDF}, we conclude that
\[
(\ev_\mbfX)_\pfwd \law \eta=\mu_\mbfX = \Dir{\ev_\mbfX(\sigma)}\comma
\]
which proves the assertion.
\end{proof}

\begin{remark}
We stress that our proof of the reverse implication in Theorem~\ref{t:MeckeDF} is different from the proof of the analogous characterization for the gamma measure, \cite[Thm.~6.3]{HagKonPasRoe13} for the case of~$\R^d$.
Indeed, the latter proof relies on a characterization of the Laplace transform of the random measure in question by some ordinary differential equation.
This approach seems however not possible in the case of the Dirichlet--Ferguson measure, the Laplace transform of which is a kind of infinite-variable hypergeometric function (see~\cite[\S4]{LzDS19}).
On the other hand, a proper analog of our proof (through the uniqueness of the solution of a multidimensional moment problem under an appropriate bound on the moments) allows one to prove the corresponding statement for the gamma measure.
\end{remark}

\begin{corollary}[Moments of the Dirichlet distribution] \label{hgr7}
Let~$\boldalpha=(\alpha_1,\dots,\alpha_k)\in\mathbb R_+^{k}$ and assume that $|\alpha|\eqdef\alpha_1+\dots+\alpha_k=1$. Then
\begin{enumerate}[$(i)$]
\item\label{i:c:Moments:1} The moments of the Dirichlet distribution $\Dir{\boldalpha}$ satisfy the following recurrence relation:
\begin{align}\label{yre7i4e}
&\int_{\R^k}\prod_{i=1}^n (\mbfs^\sym{i}\cdot\mathbf y)\,\diff\Dir{\boldalpha}(\mathbf y)\\ \notag
&\quad=\frac1n\sum_{\substack{\xi\subset\{1,\dots,n\}\\ |\xi|<n}}
\binom{n}{|\xi|}^{-1} \int_{\R^k}\prod_{i\in\xi}(\mbfs^\sym{i}\cdot\mathbf y)\,\diff\Dir{\boldalpha}(\mathbf y)\,\left(\underset{j\in\{1,\dots,n\}\setminus\xi} {\mathlarger{\mathlarger{\mathlarger\Diamond}}} \mbfs^\sym{j}\right)\cdot \boldalpha
\end{align}
for all  $n\in\mathbb N$ and $\mbfs^\sym{1},\dots,  \mbfs^\sym{n}\in\R^k$. (Here, $|\xi|$ denotes the number of elements of the set~$\xi$.) 
In particular, for all $n\in\mathbb N$ and  $\mbfs\in\R^k$,
\begin{equation}\label{trw5huy}
\int_{\R^k} (\mbfs\cdot\mathbf y)^n\,\diff\Dir{\boldalpha}(\mathbf y)=\frac1n
\sum_{i=0}^{n-1} \int_{\R^k} (\mbfs\cdot\mathbf y)^i\,\diff\Dir{\boldalpha}(\mathbf y)\big(\mbfs^{\compo (n-i)}\cdot \boldalpha\big) \fstop
\end{equation}

\item\label{i:c:Moments:2} For all $n\in\mathbb N$ and $\mbfs\in\R^k$,
\[
\int_{\R^k} (\mbfs\cdot\mathbf y)^n\,\diff\Dir{\boldalpha}(\mathbf y)= Z_n(\mbfs^{\compo 1}\cdot \boldalpha, \dotsc, \mbfs^{\compo n}\cdot \boldalpha)\comma
\]
where $Z_n$ denotes the cycle index polynomial of the symmetric group~$\mfS_n$.
\end{enumerate}
\end{corollary}

\begin{proof} Choose $\mbbX=[0,1]$, $\diff\sigma(x)=\diff x$, and choose a partition $\mbfX$ such that $\ev_\mbfX(\sigma)=\boldalpha$. Then formula \eqref{trw5huy} follows from \eqref{bgufr7rd} if we note that, for $\beta=1$,
\[
\frac{(n-1)_{n-1-i}}{(\beta+n-1)_{n-i}}=\frac{(n-1)_{n-1-i}}{(n)_{n-i}}=\frac1n\fstop
\]
 Next, note that the right hand side of formula \eqref{yre7i4e} is an $n$-linear symmetric form of $\mbfs^\sym{1},\dots,\mbfs^\sym{n}\in\R^k$, and for $\mbfs=\mbfs^\sym{1}=\dots=\mbfs^\sym{n}$, the right hand side of \eqref{yre7i4e} is equal to the right hand side of formula \eqref{trw5huy}. Hence, \eqref{yre7i4e} follows from \eqref{trw5huy} and the polarization identity. The second statement follows by noticing that the cycle index polynomials of~$\mfS_n$ satisfy the  recurrence relation \eqref{trw5huy}, e.g.~\cite[Eqn.~(2.2)]{LzDS19}.
\end{proof}

\begin{remark}\label{yr7i} Statement~\ref{i:c:Moments:2} of Corollary~\ref{hgr7} is shown in \cite[Thm.~3.3]{LzDS19} by different methods.
\end{remark}

\begin{remark}\label{4ejyr4}
By using formula \eqref{bgufr7rd}, one can immediately extend  
Statement~\ref{i:c:Moments:1} of Corollary~\ref{hgr7} to the case of a general $\boldalpha\in\mathbb R_+^{k}$.  

\end{remark}

\begin{corollary}[Moments of the Dirichlet--Ferguson measure]\label{rte6u4wu6}
Let~$\sigma\in\msP$ (i.e.~$\beta=1$) and~$\eta$ be a Dirichlet--Ferguson process. Then, the moments of~$\eta$ satisfy the following recurrence relation:
\begin{equation}\label{cyr7}
\mbbE\quadre{\prod_{i=1}^n\eta g_i } = \frac1n\sum_{\substack{\xi\subset\{1,\dots,n\}\\ |\xi|<n}}
\binom{n}{|\xi|}^{-1} \mbbE\quadre{ \prod_{i\in\xi}\eta g_i} \int \prod_{j\in\{1,\dots,n\}\setminus\xi} g_j \,\diff\sigma
\end{equation}
for all $n\in\mathbb N$ and $g_1,\dots,g_n\in \R_b(\mbbX)$. In particular, for all $n\in\mathbb N$ and  $g\in \R_b(\mbbX)$,
\begin{equation*}
\mbbE\quadre{ (\eta g)^n} =\frac1n\sum_{i=0}^{n-1}
\mbbE\quadre{ (\eta g)^{i}} \int g^{n-i}\diff\sigma \fstop
\end{equation*}
\end{corollary}

\begin{proof}In the case where the functions $g_1,\dots,g_n\in \R_b(\mbbX)$ take on a finite number of values, formula \eqref{cyr7} follows from \eqref{fst5w} and \eqref{yre7i4e}.
In the general case, formula \eqref{cyr7}  follows by approximation and the Dominated Convergence Theorem.
\end{proof}

\begin{remark}
Similarly to Remark \ref{4ejyr4}, one can easily extend Corollary \ref{rte6u4wu6} to the case of a general finite intensity measure $\sigma$.
\end{remark}

\begin{remark} A non-recursive formula for the moments of the Dirichlet--Ferguson measure, namely the full expansion of~\eqref{cyr7}, may be found in~\cite[Lemma~5.2]{Fen10}.
\end{remark}

\begin{corollary}[A characterization of the Dirichlet distribution]
Let $k\ge 2$. Let~$\theta$ be a probability measure on~$\R_+^k$. Then, the following statements are equivalent:
\begin{enumerate}[$(i)$]
\item\label{i:c:Car:1} $\theta$ is the Dirichlet distribution~$\Dir{\boldalpha}$ with parameter~$\boldalpha\in\R_+^k$;
\item\label{i:c:Car:2} for every non-negative measurable function~$g\colon \R_+^k\rar \R$,
\begin{align}\label{dye67}
\int_{\R_+^k} |\mbfy|\,g(\mbfy) \diff \theta(\mbfy) = \int_{\R_+^k} \int_0^1 (1-t)^{\abs{\boldalpha}-1} \sum_{i=1}^k \alpha_i\, g\big((1-t)\mbfy+ t\mbfe_i\big) \diff t \diff \theta(\mbfy) \fstop
\end{align}
Here $|\mbfy|\eqdef y_1+\dots+y_k$ for $\mbfy\in \R_+^k$ and $(\mbfe_i)_{i=1,\dots,k}$ is the canonical basis in~$\R^{k}$.
\end{enumerate}

Moreover, for every non-negative (or bounded) measurable function~$f\colon \Delta^{k-1}\times \{1,\dots,k\}\rar \R$,
\begin{align}\label{vhftttf}
\int_{\Delta^{k-1}} \sum_{i=1}^k y_i f(\mbfy, i) \diff \Dir{\boldalpha}(\mbfy) = \int_{\Delta^{k-1}} \int_0^1 \sum_{i=1}^k \alpha_i\, f\big((1-t) \mbfy+ t\mbfe_i, i\big) (1-t)^{\abs{\boldalpha}-1}  \diff t\diff \Dir{\boldalpha}(\mbfy) \fstop
\end{align}
\begin{proof} 
Assume~\ref{i:c:Car:1} holds. Similarly to the proof of Corollary \ref{hgr7}, choose $\mbbX=[0,1]$, $\diff\sigma(x)=\abs{\boldalpha}\diff x$, so that~$\beta=\abs{\boldalpha}$, and choose a partition $\mbfX$ such that
$\ev_\mbfX(\sigma)=\boldalpha$.
Let~$\eta$ be a Dirichlet--Ferguson process.
Applying formula \eqref{fvgyugwe} to $G\eqdef g\circ \ev_\mbfX$ 
and recalling \eqref{bur8ut} gives
\begin{align*}
&\int_{\Delta^{k-1}} g(\mbfy) \diff \Dir{\boldalpha}(\mbfy) = \mbbE\tquadre{g(\eta X_1 ,\dotsc,\eta X_k)}
\\
&=\sum_{i=1}^k \int_{X_i} \int_0^1 \mbbE\quadre{g\big((1-t)\, \eta X_1,\dots,
(1-t)\, \eta X_i+t,\dots,(1-t)\, \eta X_k\big)} (1-t)^{\abs{\boldalpha}-1} \diff t\diff x
\\
&=\int_{\Delta^{k-1}} \int_0^1 \sum_i^k \alpha_i g\big((1-t)\mbfy+ t\mbfe_i\big) (1-t)^{\abs{\boldalpha}-1} \diff t \diff \Dir{\boldalpha}(\mbfy) \fstop
\end{align*}
Thus, \eqref{dye67} holds for $\theta=\Dir{\boldalpha}$. 
Formula \eqref{vhftttf} is proven analogously by applying formula~\eqref{vcghfsky} to 
\[
F(\eta,x)=\sum_{i=1}^k f\big( \ev_\mbfX(\eta),i\big)\mathbf 1_{X_i}(x) \fstop
\]

In order to prove  that formula \eqref{dye67} uniquely identifies the measure $\theta$, one uses essentially the same arguments as in the proof of Theorem~\ref{t:MeckeDF}. One first shows that
\[
\int_{\R_+^k} \diff \theta(\mbfy) \, |\mbfy|^n=1\comma \quad n=1,2\comma
\]
which implies $|\mbfy|=1$ $\theta$-a.s., i.e., $\theta$ is concentrated on  $\Delta^{k-1}$.
Choosing $g(\mbfy)\eqdef (\mbfs\cdot \mbfy)^n$ one finds the recurrence relation for the moments of $\theta$.
\end{proof}
\end{corollary}

\section{Known applications}\label{s:Applications}
Our main result Theorem~\ref{t:MeckeDF} has already found several applications ---in particular to the study of stochastic dynamics for measure-valued processes--- some of which we briefly discuss below.

\paragraph{The Dirichlet--Ferguson diffusion}
In the case when~$\mbbX$ is a closed Riemannian manifold with volume measure~$\sigma$, the first named author constructs and studies in~\cite{LzDS22} a $\msP$-valued Markov diffusion~$\boldsymbol\eta$ with invariant measure~$\DF_\sigma$.
Theorem~\ref{t:MeckeDF} plays a crucial role in establishing an integration-by-parts formula for (whence the closability of) the Dirichlet form on~$L^2(\DF_\sigma)$ corresponding to~$\boldsymbol\eta$, thus showing that~$\boldsymbol\eta$ is the `Brownian motion' of the celebrated $L^2$-Wasserstein geometry on~$\msP$.

\paragraph{Integration by parts for the discrete gradient}
In~\cite{FliTor21}, I.~Flint and G.L.~Torrisi study a discrete gradient operator for functions on~$\msP$ defined as
\[
D_{(x,t)}F(\mu)\eqdef F\ttonde{(1-t)\mu+t\delta_x}-F(\mu)\comma \qquad x\in\mbbX\comma t\in [0,1]\comma \mu\in\msP \fstop
\]
Using our Theorem~\ref{t:MeckeDF}, the authors prove an integration-by-parts formula for the discrete gradient,~\cite[Thm.~1.1]{FliTor21}, computing the associated divergence (adjoint operator) with respect to the $L^2(\sigma\otimes \BetaD{1,\beta})$-scalar product.
This is subsequently used to obtain a Gaussian bound for the Wasserstein distance between any first-order integral of a Dirichlet--Ferguson process and a standard Gaussian random variable~\cite[Thm.~1.3]{FliTor21}, and to prove a Gaussian quantitative CLT for the first chaos~\cite[Cor.~5.3]{FliTor21}.

\paragraph{The Fleming--Viot process}
The Dirichlet--Ferguson measure~$\DF_\sigma$ is the unique stationary, reversible distribution for the Fleming--Viot process with parent-independent mutation, see, e.g.,~\cite[Thm.s~5.3,~5.4]{Fen10}.
Similarly to~\cite{LzDS22}, one can use identity~\eqref{vcghfsky} to re-establish an integration-by-parts formula for the Dirichlet form of the Fleming--Viot process which was shown by other means in~\cite[\S5.1]{OveRoeSch95}.

{\small
\bibliographystyle{plain}

}


\end{document}